# POSITIVE REALNESS OF A TRANSFER FUNCTION NEITHER IMPLIES NOR IS IMPLIED BY THE EXTERNAL POSITIVITY OF THEIR ASSOCIATE REALIZATIONS

*Abstract*—This letter discusses the differences in-between positive realness and the external positivity of dynamic systems. It is proved that both concepts do not imply mutually.

*Indexing Terms*— Dynamic systems, External positivity, Positive realness.

*Introduction*: It is commonly intended that single-input single-output linear time-invariant dynamic positive real systems (i.e. those whose transfer functions are positive real) are externally positive, [1-2], in the sense that their outputs are positive for any nonnegative inputs and zero initial conditions (see, for instance, [3-6]). However, it is proved that external positivity of a system is a more general property than its positive realness and, furthermore, some positive real systems are not externally positive. In particular, real positivity of a transfer function requires input-output stability and relative degree of at most unity while external positivity may be accomplished by certain unstable systems irrespective of the relative degree. External positive unstable systems are, for instance, those having real positive coefficients of the numerator polynomial of its transfer function while its denominator has all coefficients negative real or null except the (positive) leading one, irrespective of the relative degree. This situation leads to the most evident sufficient condition for external positivity of both continuous-time and discrete – time externally positive systems, [1]. Another essential difference is that external positivity is maintained under discretization, [1], while positive realness is never preserved under discretization of positive real continuous-time systems of unity relative degree (since positive realness of discrete-time system requires relative degree zero). Another difference, proved in this manuscript, is that inverse transfer functions of positive real transfer functions are always positive real (although they are realizable only if the original ones have zero relative degree) while inverse systems of externally positive real systems are not externally positive . Both properties are very important in Circuit Design and Systems Theory. In particular, positive realness of the uncontrolled transfer function leads to hyperstable closed –loop systems (i.e. those being stable for any feedback connection satisfying Popov´s inequality what includes the standard absolute stability problem as a particular case of asymptotic hyperstability provided that the uncontrolled open-loop transfer function is strictly positive real). The most general property of external positivity is relevant in some problems of Biology and Ecology where no variable is allowed to take negative values at any time. A common relevant property to both classes of systems is that their input-output energy is positive for all time. In this letter, the external positivity property of single-input single output continuous-time linear time-invariant dynamic systems is briefly discussed. It is found that positive realness of a transfer function neither implies nor is implied by the external positivity of its state-space realizations.

*Externally positive system* : An externally positive dynamic single- input single- output system is that satisfying the property that any input –output pair $(u(t), y(t))$ satisfies the property for zero initial conditions $u: \mathbf{R}_{0+} \to \mathbf{R}_{0+}$ ( $\mathbf{R}_{0+}$ being the set of nonnegative real numbers) , i.e. the output is non-negative for all time for any non- negative input. Equivalently, the system possesses a positive impulse response f (t) for all time.

*Positive real transfer function*: A positive (respectively, strictly positive) real transfer function $F(s)$ is that which satisfies $\operatorname{Re} F(s) \geq 0, \forall \operatorname{Re} s > 0$ (respectively, $\operatorname{Re} F(s) > 0, \forall \operatorname{Re} s \geq 0$).

Positive real transfer functions are either critically stable or stable and residuals at the critically stable poles (which are simple, if any) are nonnegative. Positive real transfer functions have absolute relative degrees of at most unity and their inverses are also positive real. Strictly positive real transfer functions are also positive real and stable (i.e. they do not possess critically stable poles) .

Consider a proper transfer function $F(s)$ and its associate impulse response $f(t)$:

$$F(s) = d + F_0(s) \; ; \; f(t) = d\delta(t) + f_0(t)$$

where d is either zero or a scalar gain and $F_0(s)$ is a strictly proper transfer function of associate impulse response $f_0(t)$, and $\delta(t)$ is the Dirac distribution. Note that the above decompositions of the transfer function and its associate impulse response are unique if d is the quotient of the leading coefficients of the numerator and denominator polynomials of the transfer function F (s) (namely, the direct input-output interconnection gain) . Note also that for $d > 0$:

$$\operatorname{Re} F(s) > 0, \forall \operatorname{Re} s \geq 0 \; \neg\Rightarrow \; \operatorname{Re} F_0(s) > 0, \forall \operatorname{Re} s \geq 0$$

where the "$\neg\Rightarrow$" symbol means " does not imply", but

$$\operatorname{Re} F(s) > 0, \forall \operatorname{Re} s \geq 0 \; \Rightarrow \; \operatorname{Re} F_0(s) > -d \text{ for all finite } \operatorname{Re} s \geq 0 \text{ since } F_0(s) \to 0 \text{ as } |s| \to \infty$$

because $F_0(s)$ is strictly proper. The following two properties are immediate from the above simple considerations:

*Property 1: A stable biproper transfer function $F(s)$ may be positive real without its associate strictly proper transfer function $F_0(s)$ being positive real.*

On the other hand,

$$f(t) = d\delta(t) + f_0(t) > 0, \forall t \geq 0$$
$$\Rightarrow d \geq 0 \land f_0(t) > 0, \forall t > 0$$

and $f_0(t) > 0, \forall t \geq 0$ if $d = 0$

*Property 2: The external positivity of a state-space realization of a (stabl, critically stable or unstable) biproper transfer function requires simultaneously a positive direct input- output interconnection gain and a positive impulse response of the associated strictly proper transfer function.*

The above properties lead to the three subsequent Assertions which prove that positive realness does not imply or it is implied by external positivity of the associate realizations and that external posiitvity is not necessarily related to stability.

*Assertion 1: External positivity does not imply positive realness*

The following example is a counterexample for Assertion 1 being false:

*Example 1*: Consider a proper transfer function $F(s) = \frac{K}{s+a}$. Its impulse response is $Ke^{-at} > 0$ for any $K > 0$ and any arbitrary real constant a for all time. In particular, for $K > 0$ and $a \geq 0$, $F(s)$ is positive real but the impulse response is positive for all finite time, and tends asymptotically to zero as time tends to infinity, even for the unstable case with a < 0 for which F(s) is not positive real. Then, the output of any state-space realization is nonnegative for zero initial conditions and any input being nonnegative for all time for any real pole and the system is externally positive as a result. Then, Assertion 1 is true.

In a more general context, it is shown in [1] that a typically inversely stable system of arbitrary order being unstable or critically stable with all the coefficients of the denominator polynomial of its transfer function being negative or null (except the leading one) is always externally positive. Assertions 2- 3 below are proved by refuting their respective converses through the same example.

*Assertion 2: Positive realness does not imply external positivity.*

In the subsequent statement, a biproper transfer function is said to be trivial if it is a pure scalar gain ( in other words if the order of any state-space realization is zero).

*Assertion 3:* Either a positive biproper real transfer function of its inverse ( which is always positive real) have not an externally positive state-space realization except in the trivial case.

The following example is a counterexample for Assertions 2-3 being false:

*Example 2:* Consider a biproper transfer function $F(s) = \frac{ds+b}{s+a} = d + \frac{b-da}{s+a}$ which is positive real for $a \geq 0, d > 0, b \geq 0$ (strictly positive real if $b > 0$ and $a > 0$). If $b \geq da$ the impulse response is positive for all time and for any nonnegative input $u(t)$, the output of any state-space realization fulfills $y(t) \geq du(t) \geq 0$, $\forall t \geq 0$ so that I the system is externally positive. However, if $b < da$, then any input $u: \mathbf{R}_{0+} \to \mathbf{R}_{0+}$ satisfying $u(t) = 0$ for some particular time t and not being identically zero on $[0,t]$ supplies an output $y(t) = -|b-da| \int_0^t f_0(t-\tau) u(\tau) d\tau < 0$ under zero initial conditions where $f_0(t)$ is the impulse response of the strictly proper transfer function $\frac{b-da}{s+a}$. Then, any state-space realization of F(s) is not externally positive while the associated transfer function is positive real. and Assertion 2 is true. The inverse of $F(s)$ is $F^{-1}(s) = \frac{s+a}{ds+b} = \frac{1}{d}\left(1 + \frac{da-b}{ds+b}\right)$ is positive real provided that F(s) is positive real and any state-space realization is externally positive if $da \geq b$ since $y(t) = d^{-1} u(t) \geq 0$, $\forall u: \mathbf{R}_{0+} \to \mathbf{R}_{0+}$, $\forall t \geq 0$. $da = b$ is the trivial case of the transfer function being a positive gain. Thus, for $da > b$ any sate-space realization associated with the positive real transfer function $F^{-1}(s)$ is externally positive while no state-space realization associated with the positive real transfer function $F(s)$ is externally positive and Assertion 3 is true.

*Acknowlwdgement:* The author is very grateful to the Spanish Ministry of Education by its financial support through Grant DPI2006-00714.